\documentclass{amsart}
\usepackage{fullpage}
\theoremstyle{plain} \newtheorem{Thm}{ }[section]
\title{Maps between moduli spaces of vector bundles and the base locus of the theta divisor}
\author{Tawanda Gwena\\ Montserrat Teixidor i Bigas}
\address{Mathematics Department, Tufts University, Medford
MA 02155, USA}

\email{tawanda.gwena@tufts.edu\\
montserrat.teixidoribigas@tufts.edu}
\begin{document}
\maketitle
\begin{section}{Introduction}

Given a vector bundle $E$ of rank $r$ and degree $d$ on a curve
$C$ of genus $g$, one can associate to $E$ in a natural way  several other
vector bundles. For example, one can take wedge powers of $E$. If
$E$ is generated by global sections, the kernel of the evaluation
map of sections is again a vector bundle. Also, new vector bundles
 can be produced by taking elementary transformations centered
at a fixed point. Under suitable conditions on degree and rank,
these constructions can be carried out globally.While all this
processes seem quite elementary, very little is known about the
resulting maps. The purpose of this paper is to fill in this gap. We
shall start by considering the kernel of the evaluation map of
sections. If $E$ is generated by global sections, consider $M_E$ the
kernel of the evaluation map of sections of $E$. Then $M_E$ can be
defined from the exact sequence
$$(*)\ \ \ \ \ \ \ 0\rightarrow M_E\rightarrow H^0(E)\otimes
{\mathcal O}_C\rightarrow E\rightarrow 0.$$

Our main result is
\begin{Thm} \label{EM_E}{\bf Theorem}
Let $d\ge (2r-1)g+1$. Then the map between the moduli spaces of
vector bundles of rank $r$ and degree $d$ and rank $d-rg$ and
degree $d$ given by $E\rightarrow M_E$ is generically injective
\end{Thm}

For the following result we assume that we are in characteristic
zero. It is known then that the wedge powers of a semistable vector
bundle are again semistable. Therefore, there is a rational map
between the moduli space $U(r,d)$ of vector bundles of rank $r$ and
degree $d$ and the moduli space of vector bundles of rank ${r\choose
i}$ and degree ${r\choose i}i{d\over r}$. We show

\begin{Thm} \label{wedge}{\bf Theorem}
Let $C$ be a projective non-singular curve defined over a field of
characteristic zero. Let $ i$ be an integer with $0<i<r$.
Then the map between the moduli space $U(r,d)$ and the moduli space
 $U\left({r\choose i},{r\choose i}i{d\over r}\right)$ given by
$E\rightarrow \wedge ^iE$ has finite fibers.

The same statement is true when replacing wedge powers by
symmetric powers.
\end{Thm}

Finally consider elementary transformations. Starting with a vector
bundle, one can produce another vector bundle of degree one more (or
one less) that essentially differs from the original one only at one
point. This construction can be done in many different ways and in
general one can get an $r$-dimensional family of such transforms
starting with a single vector bundle. We show here  that if we start
with a family of vector bundles and do all possible such transforms,
the dimension of the vector bundles so obtained never goes down.

We then give an application in finding lower bounds on the dimension
of the base locus of the theta divisor in the moduli space of vector
bundles  (see Theorem \ref{baselocus} and Corollary
\ref{baselocus2}). For many values of the rank, this improves known
results of Arcara (\cite {A}), Popa (\cite {P}) and Schneider (\cite
{S}).

\end{section}

\begin{section}{Injectivity of the map $E\rightarrow M_E$}

If $E$ is a generic vector bundle of rank $r$ and degree at least
$rg+1$, then it is generated by global sections. The same is true
for every semistable vector bundle if the degree $d$ satisfies
$d>r(2g-1)$.  Whenever $E$ is generated by global sections, one can
define the vector bundle $M_E$ as in (*). Our goal in this section
is to prove that the assignment $E\rightarrow M_E$ is generically
injective when it is globally defined (namely $d\ge 2rg+r$. What we
show in fact is a little bit more, namely that $M_{(M_E)^*}=E^*$.
This in turn is equivalent to $h^0((M_E)^*)=h^0(E)$. We shall reduce
this statement to the surjectivity of a map between spaces of
sections of vector bundles and prove the latter by induction on the
rank.

\begin{Thm}\label{noinvolucio}{\bf Lemma}
Let $C$ be a generic curve. Let $L$ be a line bundle on $C$ such
that the complete linear system $L$ has no fixed points. Then, the
map $C\rightarrow |L|$ is composed with an involution if and only
if $L=\bar L^2$ with $\bar L$ a line bundle corresponding to a
$g^1_{g+2\over 2}$ on $C$.
\end{Thm}

\begin{proof}
Assume that the map associated to the line bundle $L$ is composed
with an involution. As the curve $C$ is generic, it is not a
covering of any curve of genus at least one. Therefore the map
corresponding to the linear system $|L|$ factors through
$$C\rightarrow C' \rightarrow |L|$$
where $C'$ is a rational curve and the first map is given by a
$g^1_k$. Because of the genericity of $C$, this implies that the
corresponding Brill-Noether number $\rho= g-2(g-k+1)\ge 0$ is
positive. Equivalently, $k\ge {g+2\over 2}$. Moreover, $L=rg^1_k$
with $h^0(L)=r+1$.

 Note now that $h^1(rg^1_k)=h^0(K-rg^1_k)=0$. This is because
 $$H^0(K-rg^1_k)\subset H^0(K-2g^1_k)\subset \mbox{Ker} (
 H^0(g^1_k)\otimes H^0(K-g^1_k)\rightarrow H^0(K))$$
 and this kernel is zero by the injectivity of the Petri map for
 the generic curve.

Hence, $$r+1=h^0(L)=h^0(rg^1_k)=rk+1-g\ge r{g+2\over 2}+1-g\ge
g+r+1-g=r+1$$ Here the last inequality comes from the fact that
$r\ge 2$. Hence, all inequalities must be equalities and then
$r=2, k={g+2\over 2}$ as stated.
\end{proof}

\begin{Thm}\label{exh}
{\bf Proposition} Let $L$ be a line bundle generated by global
sections on a generic curve. Assume $L\not= \bar L^2$ where $\bar
L$ is a line bundle corresponding to a $g^1_{{g+2\over 2}}$ on
$C$. Then, the map $$\psi :H^0(L)\otimes H^0(K)\rightarrow
H^0(K\otimes L)$$ is onto. The result is also true on any curve if
$L$ is generic.
\end{Thm}
\begin{proof}
From \ref{noinvolucio}, the map corresponding to the complete
linear system $|L|$ is not composed with an involution and by
assumption it has no fixed points. The same is true on any curve
if $L$ is taken to be generic. By Riemann-Roch
$$k=h^0(L)=d+1-g+h.$$ where $h=h^1(L)=h^0(K\otimes L^{-1})$.
Choose $P_1,\ldots,P_{k-2}$ generic points on
$C$. Then,
$$h^0(L(-P_1-\cdots-P_{k-2}))=2,\ h^0(K\otimes L^{-1}(P_1+\cdots+P_{k-2}))=h^0(K\otimes L^{-1})=h.$$
 Moreover $L(-P_1-\cdots-P_{k-2})$
has no fixed points and is therefore generated by global sections.
Consider the map
$$\bar \psi :H^0(L(-P_1-\cdots-P_{k-2}))\otimes H^0(K)\rightarrow H^0(K\otimes
L(-P_1-\cdots-P_{k-2}))$$

 By the base point free pencil trick, the kernel of $\bar \psi$ is
  $H^0(K\otimes L^{-1}(P_1+\cdots+P_{k-2}))$ and by our choice this
  space has dimension $h$.
Therefore, the image of the above map has dimension $2g-h$.

Now, the $P_i$ have been chosen to  impose independent conditions
on the linear system $|L|$. Hence, they impose independent
conditions to the image of the map $\psi$. Then,
$$\dim(\mbox{Im} \psi)\ge \dim (\mbox{Im} \bar
\psi)+k-2=(2g-h)+(d+1-g+h)-2=d+g-1=h^0(K\otimes L).$$
 \end{proof}

The next result is not needed in the sequel. For the sake of
completeness, we show that the the restriction in the previous
Proposition is not arbitrary.

\begin{Thm} {\bf Lemma} Assume that $L=\bar L^2$ where $\bar L$ is a line
 bundle corresponding to a
$g^1_{{g+2\over 2}}$ on $C$. Then the morphism
$$\psi :H^0(L)\otimes H^0(K)\rightarrow H^0(K\otimes L)$$ is not onto.
\end{Thm}
\begin{proof}
  Let $s_0,s_1$ be a basis of section of $|\bar L|$. Then $s_0^2,
  s_0s_1, s_1^2$ are sections of $|L|$.

  From the base-point-free pencil trick, the map
  $$\bar \psi :<s_0^2,s_1^2>\otimes H^0(K)\rightarrow H^0(K\otimes L)$$ has
  kernel $H^0(K\otimes L^{-1})=H^0(K\otimes \bar L^{-2})$ and as
  in the proof of \ref{noinvolucio}, the latter is zero. Hence,
  the image has dimension $2g$. As $h^0(K\otimes L)=2g+1$, the map
  is not onto. In order to prove the claim, it suffices to show
  now that the image of the map $\psi$
is the same as the image of the map $\bar \psi$. This is
equivalent to proving that the elements of the form $s_0s_1t, t\in
H^0(K)$ are in the image of $\psi$.

We show first that the map
 $$<s_0,s_1>\otimes H^0(K\otimes \bar L^{-1})\rightarrow H^0(K)$$
is onto. In fact, the kernel of this map is zero because this is a
Petri map for a generic curve. Therefore the image has dimension
$$2h^0(K\otimes \bar L^{-1})=2(h^0(\bar L)+g-1-deg(\bar L))=2g$$ as
stated.

Then, given $t \in H^0(K)$, we can write
$$t=as_0+bs_1,\ a,b\in H^0(K\otimes \bar L^{-1}).$$ Then,
$$s_0s_1t=(as_1)s_0^2+(bs_0)s_1^2$$ as needed.
\end{proof}

 \begin{Thm} {\bf Theorem}. Let $C$ be a generic curve. Let $L$ be
 a line bundle on $C$ generated by global sections. Assume $L\not= \bar L^2$ where
 $\bar L$ is a line bundle corresponding to a $g^1_{{g+2\over 2}}$ on $C$.
Consider the vector bundle $M_L$ defined by the exact sequence
$$0\rightarrow M_L\rightarrow H^0(L)\otimes {\mathcal O}_C
\rightarrow L\rightarrow 0$$
 Then, $h^0((M_L)^*)=h^0(L)$. The same is true for every curve if $L$
 is taken to be generic.
\end{Thm}
\begin{proof}
Dualizing the sequence above and taking global sections, one
obtains

$$0=H^0(L^*)\rightarrow H^0(L)^*\rightarrow H^0((M_L)^*)\rightarrow
H^1(L^*)\rightarrow H^0(L)^*\otimes H^1({\mathcal O})$$ Therefore,
the statement of the proposition is equivalent to the injectivity
of the map
$$H^1(L^*)\rightarrow H^0(L)^*\otimes H^1({\mathcal O})$$
This in turn is dual of the map
$$\psi :H^0(L)\otimes H^0(K)\rightarrow H^0(K\otimes L)$$
which has been proved to be onto in \ref{exh}.
\end{proof}

\begin{Thm}\label{exhE} {\bf Proposition}
Let $E$ be a generic vector bundle of rank $r$ and degree $d\ge
(2r-1)g+1$. Then, the map
$$\psi :H^0(E)\otimes H^0(K)\rightarrow H^0(K\otimes E)$$
is onto.
\end{Thm}
\begin{proof}
By Lange's Theorem \cite{RT}, a generic vector bundle can be
written as an extension of generic stable bundles
$$0\rightarrow E_1\rightarrow E\rightarrow E_2\rightarrow 0$$
so long as the slopes satisfy

$$\mu(E_2)-\mu(E_1)\ge g-1.$$
We prove the result by induction on the rank $r$. The case $r=1$
has been proved in \ref{exh}. Assume now that $r\ge 2$ and write
$E$ as an extension
$$0\rightarrow L\rightarrow E\rightarrow E_2\rightarrow 0$$
where $L$ is a line bundle of degree $g+1$ and $E_2$ is a vector
bundle of degree at least $(2r-2)g\ge (2(r-1)-1)g+1$. Then both
$L$ and $E_2$ are generated by global sections and by induction
assumption the corresponding maps $\psi_L, \psi_{E_2}$ are onto.
Consider the exact diagram
$$\begin{matrix}0&\rightarrow &H^0(L)\otimes H^0(K)&\rightarrow &
H^0(E)\otimes H^0(K)&\rightarrow&H^0(E_2)\otimes
H^0(K)&\rightarrow &0\\
 & &\downarrow& &\downarrow & & \downarrow& & \\
 0&\rightarrow &H^0(L\otimes K)&\rightarrow &
H^0(E\otimes K)&\rightarrow&H^0(E_2\otimes K)&\rightarrow &0\\
\end{matrix}$$
The zeros on the right are obtained because $h^1(L)=0$ as $L$ is
generic of degree $g+1$ and also $h^1(K\otimes L)=0$

As the left and right vertical maps are onto, so is the middle
one.
\end{proof}

 \begin{Thm}\label{seccionsE} {\bf Theorem}. Let $C$ be a generic curve. Let $E$ be
 a generic vector bundle on $C$ of degree at least $(2r-1)g+1$.
Consider the vector bundle $M_E$ defined by the exact sequence
$$0\rightarrow M_E\rightarrow H^0(E)\otimes {\mathcal O}_C\rightarrow E\rightarrow 0$$
 Then, $h^0(M_E)=h^0(E)$.
\end{Thm}
\begin{proof}
The proof is now identical to the one in the line bundle case,
replacing \ref{exh} by \ref{exhE}.
\end{proof}



The proof of Theorem \ref{EM_E} immediately follows from the above
as $E=M^*_{(M_E)^*}$.

\end{section}

\begin{section}{Elementary transformations}
Given a vector bundle $E$ a point $P$and a linear map of the fiber
$E_P$ of $E$ at $P$ to the base field ${\bf C}$, one has an exact
sequence
$$0\rightarrow E'\rightarrow E\rightarrow {\bf C}_P\rightarrow 0.$$
Then, $E'$ is said to be obtained from $E$ by a (direct)
elementary transformation.

Dualizing the above sequence, one obtains
$$0\rightarrow E^*\rightarrow E^{'*}\rightarrow {\bf C}_P\rightarrow 0.$$
So, $E^*$ is obtained from $E^{'*}$ by a direct elementary
transformation and $E$ is said to be obtained from $E'$ by an
inverse transformation.

\begin{Thm} \label{dimtransel} {\bf Lemma}
Let $C$ be a curve, $A$ a subvariety of $U(r,d)$. Let $B$ be the set
of bundles  obtained from bundles in A by doing an elementary
transformation. Then $\dim B\ge \dim A$.
\end{Thm}
\begin{proof}
Given an $E\in A$, one obtains elements in $B$ by considering
extensions of the form
$$0\rightarrow E\rightarrow E'\rightarrow k_P\rightarrow 0$$
where $P\in C$ is a point and $k_P$ is a one dimensional
sky-scraper sheaf with support on $P$. Equivalently, one has an
exact sequence
$$0\rightarrow E^{'*}\rightarrow E^*\rightarrow k_P\rightarrow 0$$
Given $E$, one can choose a point $P$ in $C$ and a map from the
fiber of $E^*$ at $P$ to $k_P$ up to multiplication with a
constant. This then determines $E'$. Assume that $\mbox{dim}\ B<\mbox{dim}\ A$. As
every element in $A$ gives rise to an $r$-dimensional family of
elements in $B$, this would imply that every element in $B$ comes
from a family of elements in $A$ that is at least
$r+1$-dimensional. From the first exact sequence a fixed $E'$ in
$B$ comes from a family of dimension at most $r$ in $A$. Hence,
this is impossible.
\end{proof}
\end{section}

\begin{section}{Symmetric and wedge maps}

In this section, we need to work over a field of characteristic
zero, as otherwise the wedge and symmetric powers of semistable
bundles are not necessarily semistable. We prove Theorem  \ref{wedge}:

\begin{proof} The proof that follows was suggested to the authors by
S.Ramanan.

Note first that
$$\det(\wedge ^iE)=(\det E)^{\otimes {r-1 \choose i-1}}.$$
 The map
 $$\begin{matrix} \mbox{Pic}^d(C)&\rightarrow &\mbox{Pic}^{{d {r-1 \choose
 i-1}}}\\
 L&\rightarrow & L^{{\otimes {r-1 \choose i-1}}}
 \end{matrix}$$
is a finite map. Hence it suffices to prove the result for moduli
spaces of vector bundles with fixed determinant.

The Picard group of the moduli spaces of vector bundles of given
rank and determinant is ${\bf Z}$ generated by an  ample divisor
$\theta $. Consider the map
$$\begin{matrix} U(r,L)&\rightarrow &U({r \choose i}, L^{ {r-1 \choose
 i-1}})\\
 E&\rightarrow &\wedge^ {r \choose i}E
 \end{matrix}$$
 The pull-back of $\theta _{r\choose i}$ is a multiple $k\theta _r$ of the
 generator $\theta _r$ of the Picard group of $U(r,L)$ and
 therefore is ample. Hence, it intersects any positive
 dimensional subvariety. As it cannot intersect the fibers of the
 map above, it follows that the fibers are finite.

 The proof for symmetric powers is carried out in the same way.
 \end{proof}
\end{section}

\begin{section}{A lower bound for the base locus of the theta
divisor in moduli spaces of vector bundles}

Denote by $U(r, L)$ the moduli space of vector bundles of rank $r$
and fixed determinant $L$ of degree $r(g-1)$. The Picard group of
the moduli space is generated by the theta divisor that can be
described as $$\{ E\in U(r,L)|h^0(E)\ge 1\}.$$

Consider the linear system in $U(r, L)$ associated to the theta
divisor. Its base locus consists of those vector bundles $E$ such
that $h^0(E\otimes L)>0$ for every $L$ of degree zero. It was shown
first by Raynaud (see \cite{R} ) that this base locus is non-empty.
Later Popa showed in \cite{P} that it is in fact of positive
dimension. Work of Schneider and Arcara generalize Popa's results.
In all these examples, the bound on the dimension depends on the
rank (usually quite large compared with the genus). We give here a
lower bound for the dimension of the base locus which, for some
values of the rank, improves known bound.

We need some preliminary results on properties of $M_E$ that are
analogous to known properties for $M_L$ when $L$ is a line bundle.

\begin{Thm}\label{MEcontedivef} {\bf Proposition}
Let $E$ be a vector bundle of rank $r$ and degree $d$. Let $M_E$
be the dual of the kernel of the evaluation map of $E$. If $j\le
d-rg-1$ and $P_1,\ldots, P_j$ are generic points on $C$, there exists
an exact sequence
$$0\rightarrow {\mathcal O}(P_1)\oplus\cdots\oplus {\mathcal
O}(P_j)\rightarrow M_E \rightarrow F\rightarrow 0$$ where $F$ is a
vector bundle of rank $d-rg-j$.
\end{Thm}
\begin{proof} Choose generic points $P_1,\ldots,P_j$ and consider
generic one-dimensional quotients $k_{P_1},\ldots,k_{P_j}$ of
$E_{P_1},\ldots,E_{P_j}$ respectively. Let $E_1$ be the elementary
transform of $E$ associated to these quotients, namely $E_1$ is
defined by the exact sequence
$$0\rightarrow E_1\rightarrow E\rightarrow k_{P_1}
\oplus\cdots\oplus k_{P_j}\rightarrow 0.$$

Let $V$ be the kernel of the map of vector spaces
$H^0(E)\rightarrow k_{P_1} \oplus...\oplus k_{P_j}$ and
$W=H^0(E)/V$. Let $F^* $ be the kernel of the map $ V\otimes
{\mathcal O}_C\rightarrow E_1$. We get then a commutative diagram
$$\begin{matrix}0&\rightarrow &F^*&\rightarrow &V\otimes
{\mathcal O}_C&\rightarrow &E_1&\rightarrow &0\\
 & &\downarrow & &\downarrow & &\downarrow & & \\
0&\rightarrow& M_E&\rightarrow &H^0(E)\otimes {\mathcal
O}_C&\rightarrow &E&\rightarrow &0\\
 & &\downarrow & &\downarrow & &\downarrow & & \\
0&\rightarrow& A&\rightarrow& W\otimes {\mathcal O}_C&\rightarrow&
k_{P_1} \oplus\cdots\oplus k_{P_j}&\rightarrow& 0\\
\end{matrix}$$
Here $A$ denotes the cokernel of the map $F^*\rightarrow M_E$. Then,
from the lower exact row, $A={\mathcal O}(-P_1)\oplus \cdots\oplus
{\mathcal O}(-P_j)$. Dualizing the left exact sequence, the result
follows.
\end{proof}

\begin{Thm} \label{baselocus}{\bf Theorem}
Let $r$ be an  integer $r\ge 2$. Let $\beta=rg(2rg-1)$ and $L$ a
line bundle
 of degree $\beta (g-1)$.
 The theta linear system in the moduli space of vector bundles of
rank $\beta$ and fixed determinant $L$ has a base locus of
dimension at least $(r^2-1)(g-1)$.
\end{Thm}
\begin{proof}
Let $L_1$ be a line bundle of degree $3g$.  Let $E$ be a vector
bundle of rank $r$ and determinant $L_1$. Note that $E$ moves in a
family of dimension $r^2(g-1)+1-g=(r^2-1)(g-1)$. Then  $M_E$ has
rank $2rg$ and $\det((M_E)^*)=\det(E)=L_1$. Then
$$\det(\wedge ^2((M_E)^*))=(\det (M_E)^*)^{r_M-1}=
(\det (M_E)^*)^{2rg-1}=L_1^{2rg-1}.$$

From our choices, the slope of
$$\mu (\wedge ^2((M_E)^*))=2\mu ((M_E)^*)=2{3\over 2}=3.$$

 Let $P_1,\ldots,P_{g-2}, Q_1,Q_2$ be generic points on $C$. From
\ref{MEcontedivef} and the condition on $d$, ${\mathcal
O}(Q_1)\oplus {\mathcal O}(Q_2)$ is a subsheaf of $(M_E)^*$. Hence
${\mathcal O}(Q_1+Q_{2})$ is a subsheaf of $\wedge ^2 M_E$. Given a
line bundle T of degree $g-4$, there exist points
$P_1,\cdots,P_{g-2}, Q_1,Q_2$ in $C$ such that $T={\mathcal
O}(P_1+\ldots+P_{g-2}-Q_1-Q_2)$. Hence,
$$h^0(\wedge ^2(M_E)^*\otimes T)\not= 0.$$

Let $L_2$ be a line bundle such that
$L_2^{rg(2rg-1)}L_1^{2rg-1}=L$. This implies that $L_2$ has degree
$g-4$ and is determined up to a finite number of choices.

 We claim
that $\wedge ^2(M_E)^*\otimes L_2$ is in the base locus of the theta
divisor for $U(rg(2rg-1), L)$. Note first that
$$\det(\wedge
^2(M_E)^*\otimes L_2))=L_2^{\mbox{\tiny rank} (\wedge
^2(M_E)^*)}\det (\wedge ^2(M_E)^*)=L_2^{rg(2rg-1)}L_1^{2rg-1}=L.$$
Moreover, if $L_0$ is a generic line bundle of degree zero, then
$L_2\otimes L_0=T$ is a generic line bundle of degree $g-4$. Hence,
by the proof above,
$$h^0(\wedge ^2(M_E)^*\otimes L_2\otimes L_0)\not= 0$$ as needed.

Using \ref{wedge}, \ref{EM_E}, the dimension of the set of vector
bundles $E$ and the dimension of the set of vector bundles $\wedge^2
(M_E)*$ is the same, so the claim follows.
\end{proof}

\begin{Thm}\label{baselocus2}{\bf Corollary}
Let $j,\ \beta$ be defined as in the previous statement. If
$\alpha >\beta$, then the moduli space of vector bundles of rank
$\alpha$ and given determinant $L$ of degree $\alpha(g-1)$ has a
base locus of dimension at least $((r^2-1)+(\alpha -\beta )^2
)(g-1)+1$
\end{Thm}
\begin{proof} Let $E_1$ be a vector bundle of rank $\alpha -\beta$
and slope $g-1$. Let $E_2$ be a vector bundle of rank $\beta $ and
determinant $L\otimes \det(E_1)^{-1}$ in the base locus of the
theta divisor. Then, $E=E_1\oplus E_2$ is in the base locus of the
theta divisor as $h^0(E\otimes L_0)\ge h^0(E_2\otimes L_0)>0$ for
every line bundle of degree zero. Hence, our claim follows from
the previous theorem.
\end{proof}

\end{section}

\end{document}